\documentclass[11pt]{amsart}
\usepackage{lmodern}
\usepackage{amsmath, amsthm, amssymb}
\usepackage{amsfonts}
\usepackage[normalem]{ulem}
\usepackage{hyperref}

\usepackage[utf8]{inputenc}

\usepackage[margin=3.5cm]{geometry}

\usepackage{verbatim} 
\usepackage{longtable}
\usepackage{faktor}
\usepackage{mathtools}

\usepackage{tikz}
\usetikzlibrary{decorations.pathmorphing}
\tikzset{snake it/.style={decorate, decoration=snake}}

\usepackage{caption}
\usepackage{tikz-cd}
\usetikzlibrary{arrows}

\usepackage{diagbox}
\usepackage{booktabs}

\usepackage{pdflscape}

\DeclareSymbolFont{bbold}{U}{bbold}{m}{n}
\DeclareSymbolFontAlphabet{\mathbbold}{bbold}

\newtheoremstyle{mytheorem}
{3pt}
{3pt}
{\itshape}
{}
{\bf}
{.}
{.5em}
{}

\theoremstyle{mytheorem}

\newtheorem{thm}{Theorem}[section]

\theoremstyle{definition}

\theoremstyle{remark}

\DeclareFontFamily{OT1}{rsfs}{}
\DeclareFontShape{OT1}{rsfs}{n}{it}{<-> rsfs10}{}
\DeclareMathAlphabet{\curly}{OT1}{rsfs}{n}{it}

\newcommand{\Partnew}{\operatorname{Part}}

\newcommand{\ZZ}{\mathbb{Z}}

\newcommand{\rk}{\operatorname{rk}}

\newcommand{\g}{\operatorname{gcd}}

\usepackage{setspace}

\begin{document}

\title{Configurations of Points with Sum 0}
\date{\today}

\author{Christoph Schiessl}
\address{ETH Z\"urich, Department of Mathematics}
\email{christoph.schiessl@math.ethz.ch}

\thanks{ The author was supported by the grant ERC-2012-AdG-320368-MCSK.
	I want to thank Karim Adiprasato, Emanuele Delucchi, Emmanuel Kowalski, Johannes Schmitt, Junliang Shen for very helpful discussions and especially Rahul Pandharipande for his invaluable support.
	This is part of the author's PhD thesis.
}

\begin{abstract}
We compute the virtual Poincaré polynomials of the configuration space of $n$ ordered points on an elliptic curve with sum 0.
\end{abstract}

\baselineskip=14.5pt

\maketitle

\tableofcontents
For any complex quasi-projective algebraic variety $X$, the \emph{virtual Poincaré polynomial} $S(X) \in \mathbb{Z}[x]$ is defined \cite{newton}, \cite{totaroserre} by the properties 
\begin{itemize}
\item $S(X) = \sum \rk H^i(X) \, x^i$ for smooth, projective $X$,
\item $S(X) = S(X \setminus C) +S(C)$ for a closed subvariety  $C \subset X$,
\item $S(X \times Y ) = S(X) S(Y)$.
\end{itemize}

For any space $X$, the ordered configuration space 
\[ F_n(X ) = \{ x_1, \dots, x_n \in X^n | \, x_i \neq x_j \} \] is the space of $n$ distinct points in $X$. Computing its cohomology is a classical, hard problem. As the configuration space is the complement of diagonals, determining the virtual Poincaré polynomials is simpler and was done for example by Getzler \cite{getzler} \cite{getzler2}.

Let $E$ be an elliptic curve with neutral element $0$. We will compute the virtual Poincaré polynomial of the space
\[ F_n^0(E) = \{ x_1, \dots, x_n | \, x_i \neq x_j \text{ and } \sum x_i =0 \}.\]
Our approach is to decompose $F_n(X)$ in the Grothendieck ring of varieties. We use an elementary version of methods of Getzler that immediately generalizes to $F_n^0(E)$. The answer seems to be new.

The combinatorial tools are Stirling numbers and M\"obius functions and we will review them first.

\section{Stirling Numbers of First Kind}

The \emph{Stirling number of first kind} $s(n,k)$ counts the numbers of permutations in $S_n$ with exactly $k$ cycles (compare \cite[Chap. 1.3]{stanley}). Write $\Partnew( n, k)$ for all the partitions $\sigma$ of the set $\{ 1, \dots, n \}$ into $k$ disjoint, non-empty subsets $\sigma_i$. We call the $k$ subsets $\sigma_1$, \dots, $\sigma_k$ in no particular order. Then
\[ s(n,k) = \sum_{\sigma \in \Partnew(n , k)} \prod ( |\sigma_i | -1)!.\] 
Let $x$ be a positive integer. In order to determine a generating series for $s(n,k)$, we look at the action of $S_n$ on sets of functions $\{1, \dots, n \} \to \{1, \dots, x\}$. The quotient consists of the multisets of size $n$ on $\{ 1, \dots, x \}$ and has cardinality \[ \binom{n+x-1}{n} = \frac{x(x+1) \cdots (x+n-1)}{n!}.\] On the other hand, any $\tau \in S_n$ with $k$ cycles has $x^k$ fixed points. By Burnside's lemma  
\[ \frac{x(x+1) \cdots (x+n-1)}{n!}  = \frac{1}{n!} \sum_{\tau \in S_n} | \operatorname{Fix} \tau | \]
and we get 
\[ x(x+1) \cdots (x+n-1) = \sum s(n,k) x^k .\]
As it is true for all integers $x$, we have found a formal generating series.

\section{The Möbius Function of the Partition Poset}

We write $\Partnew(n)$ for the partitions of the set $\{1, \dots, n\}$. The number of parts of $\sigma \in \Partnew(n)$ is called $l(\sigma)$. The set $\Partnew(n)$ is partially ordered by setting $\sigma \le \pi$ if $\sigma$ is finer than $\pi$. Write $\mathbbold{0} = \{\{1\}, \dots, \{ n\}\}$ for the minimal partition.

\begin{thm}[M\"obius Inversion]
	For any finite poset $(M, \le)$, the M\"obius function $\mu \colon M \times M \to \ZZ$ on $M$ is defined by the relations  \begin{align*} \mu(x,z) = 0 \text{ when } x \not \le z & & \sum_{ x \le y \le z} \mu(x,y) = \delta(x,z) \text{ when } x \le z.\end{align*}
Here $\delta$ is the Kronecker delta
	\[ \delta(x,y) = \begin{cases} 1 & \text{ if } x = y \\ 0 & \text{ otherwise}. \end{cases} \]
		Let $f \colon M \to \mathbb{Z}$ a function on $M$ and 
		\[ g (x) = \sum_{x \le y} f(y).\]
Then we can reconstruct $f$ from $g$: 
	\[ f(x) = \sum_{ x \le y} \mu(x, y) g(y).\]
\end{thm}

Following \cite{BG75}, we will use M\"obius inversion to compute the M\"obius function for the poset of partitions. Let $x$ be a positive integer and $p \colon \{1, \dots, n\} \to \{ 1, \dots, x\}$ a function. The preimages of the elements of $\{ 1, \dots, x\}$ induce a partition of $ \{ 1, \dots,n \}$, that we call the kernel of $p$. Let $f(\sigma)$ be the number of functions $\{ 1, \dots, n \} \to \{1, \dots, x\}$ with kernel $\sigma$. Then $f(\mathbbold{0})$ counts all injective functions   $\{ 1, \dots, n \} \to \{1, \dots, x\}$, hence it is \[ f(\mathbbold{0}) = x(x-1) \cdots (x-n+1).\] On the other hand $g(\sigma) = \sum_{\sigma \le \pi} f(\pi)$ allows the same values on different parts on $\sigma$. Hence we have \[ g(\sigma) = x^{l(\sigma)}.\] By M\"obius inversion
\[ f(\mathbbold{0}) = \sum_{\sigma} \mu(\mathbbold{0}, \sigma) g(\sigma) \] or
\[ x(x-1) \dots (x-n+1) = \sum_{\sigma \in \Partnew(n)} \mu(\mathbbold{0}, \sigma) x^{l(\sigma)} .\]
As this holds for all values of $x$, it is valid as an identity for formal polynomials. So for the maximal partion $\mathbbold{1} = \{1, \dots, n\}$, we can immediately read off the constant term and get
\[ \mu(\mathbbold{0}, \mathbbold{1}) = (-1)^{n-1} (n-1)! .\] For general $\sigma$, the poset $\{ \pi \in \Partnew(n) | \pi \le \sigma \}$ is a product of posets
\[ \{ \pi \in \Partnew(n) | \pi \le \sigma  \} \simeq  \{ \pi \in \Partnew(|\sigma_1|) | \pi \le \sigma_1 \} \times \dots \times \{ \pi \in \Partnew(|\sigma_{l(\sigma)}|) | \pi \le \sigma_{l(\sigma)} \}\] and hence
\[ \mu(\mathbbold{0}, \sigma) = \mu (\mathbbold{0}, \sigma_1) \cdots \mu(\mathbbold{0}, \sigma_{l(\sigma)}) = (-1)^{n-l(\sigma)} \prod_i  (|\sigma_i|-1)! .\]

\section{Virtual Poincaré Polynomials of Configuration Spaces}

For any $X$, we write $[X]$ for the class of $X$ in the Grothendieck ring of varieties. 
We have maps \[F_n(X) \to F_{n-1}(X)\] with fiber $X \setminus (n-1)$ \cite{fadell}. This suggests -- ignoring possible topological problems --  \[ [F_n(X)] = [F_{n-1}(X)] \times [X - (n-1)]\] and hence \[ [F_n(X) = [X] ([X] -1 ) \cdots ([X] -n+1) = \sum_{ k \ge 0} [X]^k (-1)^{n-k} s(n,k). \] We will prove this formula be a different approach using the M\"obius function of the partition poset. It is insprired by Getzler \cite{getzler} \cite{getzler2}, who even gave a description for the $S_n$ action on 
$S(F_n(X))$.

We look at the higher diagonals \[\Delta_{\sigma} = \{x_1, \dots, x_ n \in X^n | \, x_i = x_j \text{ if $i$ and $j$ are in the same part of $\sigma$}\}\] for any partition $\sigma$ of $\{ 1, \dots, n\}$. 
By the inclusion-exclusion principle we have a decomposition \[ [F_n(X)] = [X^n] - \sum_{i \neq j} [ \{ x_i = x_j \} ] + \dots = \sum_{\sigma \in \Partnew(n)} m_\sigma [ \Delta_{\sigma}].\] for some coefficients $m_\sigma \in \mathbb{Z}$. In order to be a valid decomposition of $F_n(x)$, the coefficients $m_\sigma$ have to satisfy the condition \[\sum_{\Delta_\pi \subseteq  \Delta_\sigma} m_{\sigma} =  \begin{cases} 1 & \text{ if } \pi = 0 \\ 0 & \text{ otherwise} \end{cases} \] for any partition $\pi \in \Partnew(n)$. As $\Delta_{\pi} \subseteq \Delta_\sigma$ if and only if $\sigma \le \pi$, these equations are exactly the definition of the M\"obius function for the poset $\Partnew(n)$:
	\[ \sum_{\sigma \le \pi} \mu(\mathbbold{0},\sigma) = \begin{cases} 1 & \pi=\mathbbold{0}, \\ 0 & \text{ otherwise. } \end{cases}\]
So we get 
\[ m_{\sigma} = \mu(\mathbbold{0}, \sigma) = (-1)^{n-l(\sigma)} \prod_i ( |\sigma_i|-1)!\]
and with $[ \Delta_\sigma] = [X]^{l(\sigma)}$ we can compute:
\[ [F_n(X) ] = \sum_{\sigma \in \Partnew(n)} [X]^{l(\sigma)} (-1)^{n-l(\sigma)} \prod_i ( |\sigma_i|-1) ! = \sum_{k \ge 1} [X]^k (-1)^{n-k} s(n,k) .\]
Now applying $S$ immediately proves:
\[ S(F_n(X)) = \sum_{k \ge 1} S(X)^k (-1)^{n-k}  s(n,k) \]

\section{Configurations of Points with Sum 0}
Let $E$ be an elliptic curve with neutral element $0$. There is a map \[\Sigma \colon F_n(X) \to E, \, (x_1, \dots, x_n) \mapsto \sum x_i.\] We look at the fiber $ F_n^0(E) = \Sigma^{-1}(0) = \{ x_1, \dots, x_n \in E^n | x_i \neq x_j, \, \sum x_i = 0.\}$  

By intersecting the decomposition $[F_n(E)] = \sum_{\sigma}  m_{\sigma}  [\Delta_{\sigma}]$ with $\Sigma^{-1}(0)$ we get
\[ [F_n^0(E) ]= \sum_{\sigma \in \Partnew(n)} m_{\sigma} \, [\Delta_{\sigma} \cap \Sigma^{-1}(0)].\]
These loci have a simpler description. Take a partion $\sigma$ with $l$ parts. We see:
\begin{align*} \Delta_{\sigma} \cap \Sigma^{-1}(0)  & = \{ y_1, \dots, y_{l} \in E^{l} | \, \sum |\sigma_i| \, y_i =0\}  \end{align*}
By a coordinate change, we can compute the following solutions of this linear equation:
\begin{align*} \{ y_1, \dots, y_{l} \in E^{l} | \sum |\sigma_i| y_i = 0 \} & \simeq \{ z_1, \dots, z_{l}\in E^{l} | \gcd(|\sigma_1|, \dots, |\sigma_{l}|) z_{l} =0 \} \\ & \simeq E^{l-1} \times (\mathbb{Z}/\gcd(|\sigma_1|, \dots, |\sigma_{l}|) \mathbb{Z})^2 \end{align*}
With the notation \[ \gcd(\sigma) = \gcd(|\sigma_1|, \dots, |\sigma_{r(\sigma)}|)\] we get
\[ [F_n^0(E)] = \sum_{\sigma \in \Partnew(n)} (-1)^{n-l(\sigma)} [E]^{l(\sigma)-1} \g^2(\sigma) \prod_i (|\sigma_i|-1)! \]
Hence the following theorem is proven.
\begin{thm}
Define \[s_m(n, k) = \sum_{\sigma \in \Partnew(n, k)} \g^2(\sigma) \prod_i (|\sigma_i|-1)!. \] 
Then we have
	\[ [F_n^0(E)] = \sum_{k\ge 1} [E]^{k-1} (-1)^{n-k} s_m(n,k) \]
and 
	\[ S( F_n^0(E)) = \sum_{k \ge 1} S(E)^{k-1} (-1)^{n-k} s_m(n,k).\]
\end{thm}

The numbers $s_m(n,k)$ are a form of modified Stirling numbers. Any $\sigma \in \Partnew(n)$ with $l(\sigma) > \frac{n}{2}$  contains a part of length 1. So $\gcd(\sigma)=1$ and  \[s(n,k) = s_m(n,k)  \text{ if } k> \frac{n}{2}.\] For a prime $p$, the only partition $\sigma \in \Partnew(p)$ with $\gcd(p) \neq 1$ is $\sigma = \{\{ 1, \dots, p \}\}$. Hence \[s(p, k) = s_m(p, k) \text{ for }k>1.\] In general, \begin{align*} s(n, 1) = (n-1)! & & s_m(n, 1) = n^2 (n-1)! ,\end{align*} as $\{\{ 1, \dots, n \}\}$ is the only partition of length 1.

Unfortunately, it is not straightforward to extend the methods of \cite{getzler}, \cite{getzler2} to describe the $S_n$-action on $S(F_n^0(E))$, because the identification 
\[ \{ y_1, \dots, y_{l} \in E^{l} | \sum |\sigma_i| y_i = 0 \} \simeq E^{l-1} \times (\mathbb{Z}/\gcd(\sigma) \mathbb{Z})^2 \]
is not compatible with the $S_n$ and $S_{l}$ actions.

\newpage

\section{Tables}
Here we give the full formulas for $[F_n(E)]$ and $[F_n^0(E)]$ for all $n \le 8$.
\vspace{1cm}

\begin{tabular}{p{1cm}l}
$n$ & $[F_n(E)]$ \\ [0.5ex]
\midrule
	2 & $E^2 - E$ \\
	3 & $E^3 - 3E^2 + 2E$ \\
4 & $E^4 - 6E^3 + 11E^2 - 6E$ \\
5 & $E^5 - 10E^4 + 35E^3 - 50E^2 + 24E$ \\
6 & $E^6 - 15E^5 + 85E^4 - 225E^3 + 274E^2 - 120E$ \\
7 & $E^7 - 21E^6 + 175E^5 - 735E^4 + 1624E^3 - 1764E^2 + 720E$ \\
8 & $E^8- 28E^7 + 322E^6 - 1960E^5 + 6769E^4 - 13132E^3 + 13068E^2 - 5040E$ \\
\end{tabular}

\vspace{1cm}
\begin{tabular}{p{1cm}l}
$n$ & $[F^0_n(E)]$ \\
\midrule
2  & $E-4$ \\
3  & $E^2 - 3E + 18$ \\
4 & $E^3 - 6E^2 + 20E - 96$ \\
5  & $E^4 - 10E^3 + 35E^2 - 50E + 600$ \\
6  & $E^5 - 15E^4 + 85E^3 - 270E^2 + 864E - 4320$ \\
7  & $E^6 - 21E^5 + 175E^4 - 735E^3 + 1624E^2 - 1764E + 35280$ \\
8  & $E^7 - 28E^6 + 322E^5 - 1960E^4 + 7084E^3 - 16912E^2 + 42048E - 322560$ \\
\end{tabular}


\bibliographystyle{alpha}
\bibliography{thesis}

\end{document}